\documentclass[11pt]{amsart}

\newcommand{\w}{\omega}
\newcommand{\e}{\varepsilon}
\newcommand{\spa}{\mathrm{span}}
\newcommand{\II}{\mathbb I}
\newcommand{\diam}{\mathrm{diam}}
\newcommand{\U}{\mathcal U}
\newcommand{\id}{\mathrm{id}}
\newcommand{\conv}{\mathrm{conv}}
\newcommand{\IN}{\mathbb N}
\newcommand{\IR}{\mathbb R}
\newcommand{\IZ}{\mathbb Z}
\newcommand{\pr}{\mathrm{pr}}
\newcommand{\W}{\mathcal W}
\newcommand{\St}{\mathcal{S}t}

\newtheorem{theorem}{Theorem}
\newtheorem{lemma}{Lemma}
\newtheorem{claim}{Claim}
\newtheorem{problem}{Problem}

\title[Topological classification of closed convex sets in Fr\'echet spaces]{Topological classification of\\ closed convex sets in Fr\'echet spaces}
\author{Taras Banakh and Robert Cauty}
\address{Instytut Matematyki, Uniwersytet Humanistyczno-Przyrodniczy Jana Kochanowskiego, Kielce, Poland, and Department of Mathematics, Ivan Franko National University of Lviv, Universytetska 1, 79000, Lviv, Ukraine}
\email{t.o.banakh@gmail.com}
\address{Universit\'e Paris VI (France)}
\email{cauty@math.jussieu.fr}

\begin{document}
\begin{abstract} We prove that each non-separable completely metrizable convex subset of a Fr\'echet space is homeomorphic to a Hilbert space. This resolves an old (more than 30 years) problem of infinite-dimensional topology. Combined with the topological classification of separable convex sets due to Klee, Dobrowoslki and Toru\'nczyk, this result implies that each closed convex subset of a Fr\'echet space is homeomorphic to $[0,1]^n\times [0,1)^m\times \ell_2(\kappa)$ for some cardinals $0\le n\le\w$, $0\le m\le 1$ and $\kappa\ge 0$. 
\end{abstract}
\maketitle

The problem of topological classification of convex sets in linear metric spaces traces its history back to founders of functional analysis S.Banach and M.Fr\'echet. For separable closed  convex sets in Fr\'echet  spaces this problem was resolved by combined efforts of V.Klee \cite{Klee} (see \cite[III.7.1]{BP}), Dobrowolski and Toru\'nczyk \cite{DT1}, \cite{DT2}:

\begin{theorem}[Klee-Dobrowolski-Toru\'nczyk]\label{KDT} Each separable closed convex subset $C$ of a Fr\'echet space is homeomorphic to $[0,1]^n\times [0,1)^m\times(0,1)^k$ for some cardinals $0\le n,k\le\w$ and $0\le m\le 1$. In particular, $C$ is homeomorphic to the separable Hilbert space $l_2$ if and only if $C$ is not locally compact.
\end{theorem}

By a {\em Fr\'echet space} we mean a locally convex complete linear metric space. A {\em linear metric space} is a linear topological space endowed with an invariant metric that generates its topology. A topological space is called {\em completely metrizable} if its topology is generated by a complete metric.

In this paper we study the topological structure of non-separable (completely metrizable) convex sets in Fr\'echet spaces and prove the following theorem that  answers an old problem LS10 posed in Geoghegan's list \cite{Ge} and then repeated in \cite{We} and \cite{BZ}.

\begin{theorem}\label{t1} Each non-separable completely metrizable convex subset of a Fr\'echet space is homeomorphic to a Hilbert space.
\end{theorem} 

Theorems~\ref{KDT} and \ref{t1} imply the following topological classification of closed convex subset in Fr\'echet spaces.

\begin{theorem} Each closed convex subset $C$ of a Fr\'echet space is homeomorphic to $[0,1]^n\times[0,1)^m\times \ell_2(\kappa)$ for some cardinals $0\le n\le\w$, $0\le m\le 1$ and  $\kappa\ge 0$. In particular, $C$ is homeomorphic to an infinite-dimensional Hilbert space if and only if $C$ is not locally compact.
\end{theorem}

Here $\ell_2(\kappa)$ stands for the Hilbert space that has an orthonormal basis of cardinality $\kappa$. 
The topology of infinite-dimensional Hilbert spaces was characterized by Toru\'nczyk \cite{Tor1}, \cite{Tor2}. This characterization was used in the proof of the following criterion from \cite{BZ} which is our main tool for the proof of Theorem~\ref{t1}.

\begin{theorem}[Banakh-Zarichnyy]\label{BZ} A convex subset $C$ in a linear metric space is homeomorphic to an infinite-dimensional Hilbert space if and only if $C$ is a completely metrizable absolute retract with LFAP.
\end{theorem}

A topological space $X$ is defined to have the {\em locally finite approximation property} (briefly, LFAP), if for each open cover $\U$ of $X$ there is a sequence of maps $f_n:X\to X$, $n\in\w$, such that each $f_n$ is $\U$-near to the identity $\id_X:X\to X$ and the family $\big(f_n(X)\big)_{n\in\w}$ is locally finite in $X$. The latter means that each point $x\in X$ has a neighborhood $O(x)\subset X$ that meets only finitely many sets $f_n(X)$, $n\in\w$. 

Theorem~\ref{t1} follows immediately from Theorem~\ref{BZ}, the Borsuk-Dugundji Theorem \cite[II.3.1]{BP} (saying that convex subsets of Fr\'echet spaces are absolute retracts) and the following theorem that will be proved in Section~\ref{s:pf}. 

\begin{theorem}\label{t2} Each non-separable convex subset of a Fr\'echet space  has LFAP.
\end{theorem}

\section{Separated Approximation Property}

Theorem~\ref{t2} establishing LFAP in non-separable convex sets will be proved with help of the metric counterpart of LFAP, called SAP.

A metric space $(X,d)$ is defined to have the {\em separated approximation property} (briefly, SAP) if for each $\e>0$ there is a sequence of maps $f_n:X\to X$, $n\in\w$, such that each $f_n$ is $\e$-homotopic to $\id_X$ and the family $(f_n(X))_{n\in\w}$ is {\em separated} in the sense that $\inf_{n\ne m}d(f_n(X),f_m(X))>0$. 

Here for two non-empty subsets $A,B\subset X$ we put $d(A,B)=\inf\{d(a,b):a\in A,\;b\in B\}$. Two maps $f,g:A\to X$ are called {\em $\e$-homotopic} if they can be linked by a homotopy $(h_t)_{t\in\II}:A\to X$ such that $h_0=f$, $h_1=g$ and $\diam \{h_t(a):t\in\II\}\le \e$ for all $a\in A$. By $\II$ we denote the unit interval $[0,1]$.

The following lemma is proved by analogy with Lemma 1 of \cite{DT2} and Lemma~5.2 of \cite{BZ}.

\begin{lemma}\label{SAP-LFAP} Each metric space with SAP satisfies LFAP.
\end{lemma}

\begin{proof} Assume that a metric space $(X,d)$ has SAP. To show that $X$ has LFAP, fix an open cover $\U$ of $X$ and find a non-expanding function $\e:X\to(0,1)$ such that the cover $\{B_d(x,\e(x)):x\in X\}$ refines the cover $\U$. For every $k\in\w$ consider the closed subset $X_k=\{x\in X:\e(x)\ge 2^{-k}\}$ of  $X$. Put $\e_k=1/4^{k+2}$ for $k\le 1$ and let $f_0:X\times\w\to X$, $f_0:(x,n)\mapsto x$, be the projection. By induction we shall construct a sequence $(\e_k)_{k\in\w}$ of positive real numbers and a sequence of maps $f_k:X\times\w\to X$, $k\in\w$, such that
the following conditions are satisfied:
\begin{enumerate}
\item $\e_k\le \frac14 \e_{k-1}\le \frac1{4^{k+2}}$;
\item $f_{k}$ is $\e_k$-homotopic to $f_{k-1}$;
\item $f_{k}|X_{k-3}\times\w=f_{k-1}|X_{k-3}\times\w$;
\item $f_k|(X\setminus X_{k+1})\times\w=f_0|(X\setminus X_{k+1})\times\w$;
\item $\inf_{n\ne m}d(f_{k}(X_{k}\times\{n\}),f_{k}(X_{k}\times\{m\}))\ge4\e_{k+1}$.
\end{enumerate}

Assume that maps $f_{i}:X\times\w \to X$ and numbers $\e_{i+1}$ satisfying the conditions (1)--(4) have been constructed for all $i<k$. By SAP, there is an $\e_k$-homotopy $(h_t)_{t\in\II}:X\times\w\to X$ such that $h_0=f_0$ and  $\delta=\inf_{n\ne m}d(h_1(X\times\{n\}),h_1(X\times\{m\})>0$. Choose a continuous function $\lambda:X\to[0,1]$ such that $X_{k-3}\cup(X\setminus X_{k+1})\subset \lambda^{-1}(0)$ and $X_k\setminus X_{k-2}\subset\lambda^{-1}(1)$.
Take any positive number $\e_{k+1}\le\frac14\min\{\delta,\e_k\}$ and define a function $f_k:X\times\w\to X$ by 
$$f_k(x,n)=h_{\lambda(x)}(f_{k-1}(x,n),n).$$
It is clear that the conditions (1)--(4) are satisfied. The condition (5) will follow as soon as we check that $d(f_k(x,n),f_k(y,m))\ge 4\e_{k+1}$ for any points $x,y\in X_k$ and distinct numbers $n\ne m$.

Find unique numbers $i,j\le k$ such that $x\in X_i\setminus X_{i-1}$ and $y\in X_j\setminus Y_{j-1}$. If $i,j<k$, then $$d(f_k(x,n),f_k(y,m))\ge d(f_{k-1}(x,n),f_{k-1}(y,m))-2\e_k\ge 4\e_k-2\e_k=2\e_k\ge 4\e_{k+1}.$$

It remains to consider the case $\max\{i,j\}=k$. We lose no generality assuming that $i=k$. If $j\ge k-1$, then 
$$d(f_k(x,n),f_k(y,m))=d(h_1(f_{k-1}(x,n),n),h_1(f_{k-1}(y,m),m))\ge\delta\ge 4\e_{k+1}.$$
Next, assume that $j\le k-2$. In this case $k\ge j+2\ge 3$. Then $$\e(x)<2^{-i+1}=2^{-k+1}<2^{-k+2}\le 2^{-j}\le\e(y)$$ and the non-expanding property of $\e$ imply that $d(x,y)\ge |\e(x)-\e(y)|\ge 2^{-j}-2^{-k+1}\ge 2^{-j-1}.$
It follows from (4) and (2) that $$d(x,f_k(x,n))=d(f_{i-2}(x,n),f_k(x,n))=d(f_{k-2}(x,n),f_k(x,n))\le \e_{k-1}+\e_k\le 2\e_{k-1}\le\frac{2}{4^{k+1}}$$and
$$d(y,f_k(y,m))=d(f_{j-2}(y,m),f_k(y,m))\le \e_{j-1}+\dots+\e_k\le 2\e_{j-1}\le\frac2{4^{j+1}}.$$
Then
$$
\begin{aligned}
d(f_k(x,n),f_k(y,m))&\ge d(x,y)-d(x,f_k(x,n))-d(y,f_k(y,m))\\
&\ge \frac1{2^{j+1}}-\frac2{4^{k+1}}-\frac2{4^{j+1}}\ge
\frac1{2^{j+1}}-\frac2{4^{j+3}}-\frac2{4^{j+1}}\ge \frac4{4^{j+5}}\ge 
 \frac4{4^{k+3}}\ge 4\e_{k+1}.
\end{aligned}
$$
This completes the inductive step.
\smallskip

After completing the inductive construction, let $f_\infty=\lim_{k\to\infty}f_k:X\times\w\to X$. The conditions (1)-(3) guarantee that the limit function $f_\infty$ is well-defined and continuous.
Let us show that $f_\infty$ is $\e$-near to $f_0$. Given any point $(x,n)\in X\times \w$, find a unique number $i\in\IN$ such that $x\in X_i\setminus X_{i-1}$. By (3) and (4), $f_\infty(x,n)=f_{i+2}(x,n)$ and $f_0(x,n)=f_{i-2}(x,n)$. Then 
$$
\begin{aligned}
d(f_\infty(x,n),x)&=d(f_\infty(x,n),f_0(x,n))=d(f_{i+2}(x,n),f_{i-2}(x,n))\\
&\le e_{i+2}+\dots+\e_{i-1}\le 2\e_{i-1}\le\frac2{4^i}<\frac1{2^i}\le \e(x).
\end{aligned}$$
The choice of the function $\e$ guarantees that $f_\infty$ is $\U$-near to the projection $f_0:X\times\w\to X$. 

It remains to prove that the family $(f_\infty(X\times\{n\})_{n\in\w}$ is discrete in $X$. Given any point $x\in X$, let $i\in\IN$ be the unique number such that $x\in X_i\setminus X_{i-1}$. Consider the ball $B(x;1/2^{i+2})=\{x'\in X:d(x,x')<1/2^{i+2}\}$ centered at $x$. 

\begin{claim} $B(x;1/2^{i+2})\cap f_\infty(X\times\w)\subset f_\infty(X_{i+1}\times\w)$.
\end{claim}
 
\begin{proof} Assume conversely that $f_\infty(y,m)\in O(x)$ for some $y\in X\setminus X_{i+1}$ and $m\in\w$. Let $j\in\w$ be a unique number with $y\in X_j\setminus X_{j-1}$. It follows from $y\notin X_{i+1}$ that $j\ge i+2$. 
Since $$d(f_\infty(y,m),y)=d(f_{j+3}(y,m),f_{j-2}(y,m))\le 2\e_{j-1}\le\frac2{4^{j+1}}\le \frac1{2^{i+2}},$$ and
$\e(y)<\frac1{2^{j-1}}<\frac1{2^i}\le \e(x)$, by the non-expanding property of $\e$, we get a contradiction:
$$\frac1{2^{i+1}}\le\frac1{2^i}-\frac1{2^{j-1}}\le|\e(x)-\e(y)|\le d(x,y)\le d(x,f_\infty(y,m))+d(f_\infty(y,m),y)<\frac1{2^{i+2}}+\frac1{2^{i+2}}=\frac1{2^{i+1}}.$$  
\end{proof}

Now the condition (5) and the inequality $\e_{i+2}\le\frac1{4^{i+4}}\le\frac1{2^{i+2}}$ implies that the ball $B(x;\e_{i+2})$ meets at most one set $f_\infty(X_{i+1}\times \{n\})$ and hence at most one set $f_\infty(X\times\{n\})$, which means that the family $\big(f_\infty(X\times\{n\})\big)_{n\in\w}$ is discrete in $X$ and hence $X$ has LFAP.
\end{proof}

\section{SAP in non-separable convex cones}

In this section we shall prove that non-separable convex cones in Fr\'echet spaces have SAP.

A subset $C$ of a linear metric space $(L,d)$ is called a {\em convex cone} if it is convex and $\IR_+\cdot C=C$ where $\IR_+=[0,\infty)$. 
The principal result of this section is

\begin{lemma}\label{coneSAP} Each non-separable convex cone $C$ in a Fr\'echet space $L$ has SAP.
\end{lemma}

For the proof of this lemma we shall use an operator version of Josefson-Nissenzweig Theorem proved in \cite{BB}:

\begin{lemma}\label{BB} For any dense continuous non-compact linear operator $S:X\to Y$ between normed spaces there is a linear continuous operator $T:Y\to c_0$ such that the operator $TS:X\to c_0$ is not compact.
\end{lemma}  
  
Let us recall that an operator $T:X\to Y$ between linear topological spaces is
\begin{itemize}
\item {\em dense} if $TX$ is dense in $Y$;
\item {\em compact} if the image $T(U)$ of some open neighborhood $U\subset X$ of zero is totally bounded in $Y$. 
\end{itemize}
A subset $B$ of a linear topological space $Y$ is {\em totally bounded} if for each open neighborhood $V\subset Y$ of zero there is a finite subset $F\subset Y$ such that $B\subset V+F$. 
\medskip

\noindent{\em Proof of Lemma~\ref{coneSAP}.} Assume that $C$ is a non-separable convex cone in a Fr\'echet space $L$. By \cite[I.6.4]{BP}, the topology of the Fr\'echet space $L$ is generated by an invariant metric $d_L$ such that for every $\e>0$ the $\e$-ball $B_L(\e)=\{x\in L:d_L(x,0)<\e\}$ centered at the origin is convex. We lose no generality assuming that the linear subspace $C-C$ is dense in the Fr\'echet space $L$.

Given any $\e>0$,  we need to construct maps $f_k:C\to C$, $k\in\w$, such that each $f_k$ is $\e$-homotopic to $\id_C$ and $\inf_{k\ne n}d(f_k(C),f_n(C))>0$.  
Since the metric $d$ has convex balls, any two $\e$-near maps into $C$ are $\e$-homotopic.

\begin{claim} There is a linear continuous operator $R:L\to Y$ onto a normed space $Y$ such that the image $R(C)$ is not separable and  
$R^{-1}(\bar B_Y)\subset B_L(\e/2)$ where $\bar B_Y=\{y\in Y:\|y\|\le 1\}$ is the closed unit ball in the normed space $Y$.
\end{claim}  

\begin{proof} By \cite[I.6.4]{BP}, the Fr\'echet space $L$ can be identified with a closed linear subspace of the countable product $\prod_{i\in\w}X_i$ of Banach spaces. For every $n\in\w$ let $Y_n=\prod_{i<n}X_i$ and $\pr_n:L\to X_n$ be the natural projection. Since $C$ is non-separable, there is $n\in\w$ such that for every $m\ge n$ the image 
$\pr_m(C)\subset Y_m$ is not separable. We can take $m\ge n$ so large that
the ball $B_L(\e/2)$ contains the preimage $\pr_m^{-1}(U)$ of some open neighborhood $U\subset Y_m$ of the origin. The neighborhood $U$ contains the closed $r$-ball $\bar B_{Y_m}(r)$ of the Banach space $Y_m$ for some $r>0$. Finally, consider the linear space $Y=\pr_m(L)\subset Y_m$ endowed with the norm $\|y\|=\frac1r\|y\|_m$ where $\|\cdot\|_m$ is the norm of the Banach space $Y_m$. Then the operator $R=\pr_m:L\to Y$ has the desired properties.
\end{proof}

In the convex cone $C$ consider the convex subset $B_C=C\cap R^{-1}(\bar B_Y)$ and observe that $C=\IR_+\cdot C_1$ and hence $T(C)=\IR_+\cdot T(B_C)$. Since the space $T(C)$ is not separable, $T(B_C)$ is not separable too. Consider the convex bounded symmetric subset $D=T(B_C)-T(B_C)\subset Y$ and observe that $\IR\cdot D=R(C)-R(C)=R(C-C)$. Then the Minkowski functional 
$$\|x\|_Z=\inf\{\lambda>0:x\in \lambda D\}$$is a well-defined norm on the linear space $Z=\IR\cdot D=R(C-C)$ and the identity inclusion $I:Z\to Y$ is a bounded linear operator from the normed space $(Z,\|\cdot\|_Z)$ to the Banach space $Y$. Since $I(Z)=Z$ is non-separable, the operator $I$ is not compact. By Lemma~\ref{BB}, there is a bounded operator $T:Y\to c_0$ such that the composition $TI:Z\to c_0$ is not compact. The latter means that the image $T(D)=TR(B_C)-TR(B_C)$ is not totally bounded in $c_0$ and hence the bounded set $TR(B_C)$ is not totally bounded in $c_0$. 

Consequently, there is $\delta\in(0,1]$ such that for every $n\in\w$ 
\begin{equation}\label{eq2}
TR(B_C)\not\subset\{(x_i)_{i\in\w}\in c_0:\max_{i\ge n}|x_i|<\delta\}.
\end{equation}
For every $n\in\w$ let $e_n^*\in c_0^*$, $e_n^*:(x_i)_{i\in\w}\mapsto x_n$, be the $n$th coordinate functional of $c_0$ and let $z^*_n=(TR)^*(e^*_n)\in L^*$. 

\begin{claim}\label{cl3} There are an increasing number sequence $(m_k)_{k\in\w}$ and a sequence $(z_k)_{k\in\w}\subset B_C$ such that for every $k\in\w$;
\begin{enumerate}
\item $|z^*_{m_k}(z_k)|\ge\delta$;
\item $|z^*_{m_i}(z_k)|<\delta^3/100$ for all $i>k$.
\end{enumerate}
\end{claim}

\begin{proof} The sequences $(m_k)$ and $(z_k)$ will be constructed by induction. By (\ref{eq2}) there are a point $z_0\in B_C$ and a number $m_0\in\w$ such that $|e^*_{m_0}(z_0)|\ge\delta$. Now assume that for some $k\in\w$ points $z_0,\dots,z_k$ and numbers $m_0<m_1<\cdots<m_k$ have been constructed. Since the points $TR(z_i)$, $i\le k$, belong to the Banach space $c_0$, there is a number $m>m_k$ so large that $|e^*_n(TR(z_i))|<\delta^3/100$ for all $n\ge m$ and $i\le k$. By (\ref{eq2}), there are a point $z_{k+1}\in B_C$ and a number $m_{k+1}\ge m$ such that $|z^*_{m_k}(z_{k+1})|=|e^*_{m_k}(TR(z_{k+1}))|\ge\delta$. This complete the inductive step.
\end{proof}  

Divide $\w$ into the countable union $\w=\bigcup_{k\in\w}N_k$ of  pairwise disjoint infinite subsets and by induction define a function $\xi:\w\times\w\to\w$ such that
$\xi(i,k)\in N_k$ and $\xi(i+1,k)>\xi(i)> i$ for all $i,k\in\w$. For any numbers $i,k\in\w$ let $$z_{i,k}:=z_{\xi(i,k)}\mbox{ and }z_{i,k}^*:=z^*_{m_{\xi(i,k)}}=(TR)^*(e^*_{m_{\xi(i,k)}}),$$ where $(z_i)_{i\in\w}$  and $(m_k)_{k\in\w}$ are given by Claim~\ref{cl3}.
It follows that the double sequences $(z_{i,k})_{i,k\in\w}$ and 
$(z^*_{i,k})_{i,k\in\w}$ have the following properties (that will be used in the proof of Claim~\ref{cl8} below):

\begin{claim} If $(i,k),(j,n)\in\w\times\w$, then
\begin{enumerate}
\item $|z^*_{i,k}(z_{i,k})|\ge\delta$;
\item $|z^*_{j,k}(z_{i,n})|<\delta^3/100$ provided $\xi(j,k)>\xi(i,n)$;
\item $|z^*_{i,k}(z)|\le 1$ for any $z\in B_C$.
\end{enumerate}
\end{claim}

\begin{claim}\label{cl4} There is a map $f:C\to C$ such that $d(f,\id)<\e/2$ and each point $x\in C$ has a neighborhood $O(x)$ whose image $f(O(x))$ lies in the convex hull $\conv(F_x)$ of some finite subset $F_x\subset C$.
\end{claim}

\begin{proof} Using the paracompactness of the metrizable space $C$, find a locally finite open cover $\U$ of $X$ that refines the cover of $C$ by open $\frac{\e}4$-balls. In each set $U\in\U$ pick up a point $c_U\in U$. Let $\{\lambda_U:C\to[0,1]\}_{U\in\U}$ be a partition of the unity, subordinated to the cover $\U$ in the sense that $\lambda_U^{-1}((0,1])\subset U$ for all $U\in\U$. Finally, define a map $f:C\to C$ by the formula $$f(x)=\sum_{U\in\U}\lambda_U(x)c_U.$$
It is standard to check that $f$ has the desired property.
\end{proof}

For every $k\in\IZ$ by $C_k$ denote the set of points $x\in C$ that have neighborhood $O(x)\subset C$ such that for each point $x'\in O(x)$ and a non-negative number $m\ge k$ we get $|z^*_mf(x')|<\delta^3/100$. It is clear that each set $C_k$ is open in $C$ and lies in $C_{k+1}$.

\begin{claim} $C=\bigcup_{k\in\w}C_k$.
\end{claim}

\begin{proof} By Claim~\ref{cl4}, each point $x\in C$ has a neighborhood $O(x)\subset C$ such that $f(O(x))\subset\conv(F)$ for some finite subset $F\subset C$. Taking into account that $TR(F)$ is a finite subset of the Banach space $c_0$,  we can find a number $m\in\w$ such that $|e^*_nTR(z)|<\delta^3/100$ for all $n\ge m$ and all $z\in F$. Then also $|e^*_nTR(z)|<\delta^3/100$ for all $z\in\conv(F)$, in particular, $|e^*_nTRf(x')|<\delta^3/100$ for any $x'\in O(x)$. This means that $x\in C_m$ by the definition of the set $C_m$.
\end{proof}

\begin{claim} There is an open cover $(U_k)_{k\in\w}$ of the space $C$ such that  $U_k\subset \bar U_k\subset C_{k-1}\cap U_{k+1}$ for all $k\in\w$.
\end{claim}

\begin{proof} By Theorem 5.2.3 of \cite{En}, there is an open cover $(V_k)_{k\in\w}$ of $X$ such that $\bar V_k\subset C_k\cap V_{k+1}$ for all $k\in\w$. For each $x\in X$ find the smallest number $k\in\w$ with $x\in C_k$  and the largest number $n\le k$ with $x\notin \bar V_n$ and put $O(x)=C_k\setminus \bar V_n$. 
Consider the open cover $\W_0=\{O(x):x\in C\}$ and observe that $\St(\bar V_k,\W_0)\subset C_k$ for every $k\in\w$. Here $\St(A,\W_0)=\cup\{W\in\W_0:W\cap A\ne\emptyset\}$ for a subset $A\subset C$.

Using the paracompactness of the space $C$, for every $n\in\w$ by induction find an open cover $\W_n$ of $C$ whose star $\St(\W_{n+1})=\{\St(W,\W_{n+1}):W\in\W_{n+1}\}$ is inscribed into the cover $\W_n$. Then the open sets $$U_k=\St(\bar V_{k-1}\cup U_{k-1},\W_{k+1}),\;\;k\in\w,$$ have the required property:
$U_k\subset\bar U_k\subset C_{k-1}\cap U_{k+1}$ for all $k\in\w$.
\end{proof}

By Theorem 5.1.9 of \cite{En} there us a partition of unity $\{\lambda_k:C\to[0,1])_{k\in\w}$, subordinated to the cover $\{U_{k+1}\setminus\bar U_{k-1}\}_{k\in\w}$ of $C$ in the sense that $\lambda_k^{-1}(0,1]\subset U_{k+1}\setminus\bar U_{k-1}$ for all $k\in\w$ (here we assume that $U_k=\emptyset$ for $k<0$). 
  
Now, for every $k\in\w$ define a map
$f_k:C\to C$ by the formula
$$f_k(x)=f(x)+\sum_{i\in\w}\lambda_i(x)z_{i,k}=f(x)+\lambda_i(x)z_{i,k}+
(1-\lambda_i(x))z_{i+1,k},$$where $i$ is the unique number such that $x\in U_{i+1}\setminus U_i$.
Since $f_k(x)-f(x)\in B_C\subset B_d(\e/2)$, we conclude that $d(f(x),f_k(x))<\e/2$ and hence $$d(x,f_k(x))\le d(x,f(x))+d(f(x),f_k(x))<\frac{\e}2+\frac{\e}2=\e$$ for all $x\in C$. So, each function $f_k:C\to C$ is $\e$-near and $\e$-homotopic to the identity $\id_C:C\to C$.

\begin{claim}\label{cl8} The family $(f_k(C))_{k\in\w}$ is separated.
\end{claim}

\begin{proof}  By the continuity of the operator $TR:L\to c_0$, there is $\eta>0$ such that $TR(B_L(\eta))\subset B_{c_0}(\delta^3/20)$. We claim that $\inf_{n\ne k}d(f_n(C),f_k(C))\ge\eta$.

Fix any distinct numbers $n,k\in\w$ and points $x,y\in C$. By the choice of $\eta$, the inequality $d(f_k(x),f_n(y))\ge\eta$ will follow as soon as we check that $\|TR(f_k(x)-f_n(y))\|>\delta^3/20$. The latter inequality will follow as soon as we find $m\in\w$ such that 
$|e^*_mTR(f_k(x)-f_n(y))|>\delta^3/20$. Since $e^*_mTR(z)=z^*_m(z)$ for all $z\in L$, it suffices to show that $|z^*_m(f_k(x)-f_n(y))|>\delta^3/20$ 
for some $m\in\w$. 

 Since $C=\bigcup_{i\in\w}U_{i+1}\setminus U_i$, there are unique numbers $i,j\in\w$ such that $x\in U_{i+1}\setminus U_i$ and $y\in U_{j+1}\setminus U_j$.
Then 
$$
\begin{aligned}
f_k(x)=f(x)+\lambda_i(x)z_{i,k}+\lambda_{i+1}(x)z_{i+1,k},\\
f_n(y)=f(y)+\lambda_j(y)z_{j,n}+\lambda_{j+1}(y)z_{j+1,n}.
\end{aligned}
$$
Without loss of generality, $\xi(i+1,k)<\xi(j+1,n)$. 

Since $x,y\in U_{\max\{i,j\}+1}\subset C_{\max\{i,j\}}$, we conclude that
\begin{equation}\label{eq1}
\max\{|z^*_m(f(x))|,|z^*_m(f(y))|\}<\frac{\delta^3}{100}\mbox{ for all }m\ge\max\{i,j\}
\end{equation}according to the definition of the set $C_{\max\{i,j\}}$.

We shall consider five cases.

1) $\lambda_{j+1}(y)>\delta^2/10$. In this case, put $m=m_{\xi(j+1,n)}$ and observe that $|z^*_m(z_{j+1,n})|=|z^*_{j+1,n}(z_{j+1,n})|\ge\delta$. Since $\max\{\xi(j,n),\xi(i+1,k),\xi(i,k)\}<\xi(j+1,k)$, we conclude that 
$$\max\{|z^*_m(z_{j,n})|,|z^*_m(z_{i+1,k})|,|z^*_m(z_{i,k})|\}<\delta^3/100$$by Claim~\ref{cl4}.
It follows from (\ref{eq1}) and $\max\{i,j\}\le\max\{\xi(i,n),\xi(j,k)\}$ that $$\max\{|z^*_m(f(x))|,|z^*_m(f(y))|\}<\frac{\delta^3}{100}.$$
Now we see that
$$
\begin{aligned}
|z^*_m(f_n(y)-f_k(x))|&=
|z^*_m(\lambda_{j+1}(y)z_{j+1,n}+\lambda_j(y)z_{j,n}+f(y)-f(x)-
\lambda_{i+1}(x)z_{i+1,k}-\lambda_i(x)z_{i,k})|\\
&\ge \lambda_{j+1}(y)|z^*_m(z_{j+1,n})|-|z^*_m(\lambda_j(y)z_{j,n}+f(y)-f(x)-
\lambda_i(x)z_{i,k}+\lambda_{i+1}(x))z_{i+1,k})|\\
&> \frac{\delta^2}{10}\delta-5\frac{\delta^3}{100}\ge \frac{\delta^3}{20}.
\end{aligned}
$$

2) $\lambda_{j+1}(y)\le \delta^2/10$ and $\xi(j,n)>\xi(i+1,k)$. In this case put  $m=\xi(j,n)$. Arguing as in the preceding case, we can show that $\max\{|z^*_m(f(x))|,|z^*_m(f(y))|\}<\delta^3/100$ and
$\max\{z^*_m(z_{i+1,k})|,|z^*_m(z_{i,k}|\}<\delta^3/100$. 
Then
$$
\begin{aligned}
|e^*_m(f_n(y)-f_k(x))|&=|e^*_m(\lambda_{j}(y)z_{j,n}+\lambda_{j+1}(y)z_{j+1,n}+f(y)-f(x)-\lambda_{i+1}(x)z_{i+1,k}-\lambda_i(x)z_{i,k})|\\
&\ge \lambda_{j}(y)|z^*_m(z_{j,n})|-\lambda_{j+1}(y)|z^*_m(z_{j+1,n})|-
|z^*_m(f(y)-f(x)-
\lambda_i(x)z_{i,k}+\lambda_{i+1}(x)z_{i+1,k})|\\
&\ge (1-\lambda_{j+1}(y))\delta-\frac{\delta^2}{10}-4\frac{\delta^3}{100}\ge
(1-\frac{\delta^2}{10})\delta-\frac{\delta^2}{10}-\frac{\delta^3}{25}>\frac{\delta^3}{20}.
\end{aligned}
$$

3)  $\lambda_{j+1}(y)\le \delta^2/10$, $\xi(j,n)<\xi(i+1,k)$, and $\lambda_{i+1}(x)>\delta/4$. In this case put $m=\xi(i+1,n)$ and observe that 
$$
\begin{aligned}
|z^*_m(f_k(x)-f_n(y))|&\ge \lambda_{i+1}(x)|z^*_m(z_{i+1,k})|-\lambda_{j+1}(y)|z^*_m(z_{j+1,n})|-|z^*_m(f(x)+\lambda_i(x)z_{i,k}-f(y)-\lambda_j(y)z_{j,n})|\\
&>\frac{\delta}4\delta-\frac{\delta^2}{10}-4\frac{\delta^3}{100}>\frac{\delta^3}{20}.
\end{aligned}
$$

4) $\lambda_{j+1}(y)\le\delta^2/10$, $\xi(j,n)<\xi(i+1,k)$, $\lambda_{i+1}(x)\le \delta/4$, and $\xi(i,k)<\xi(j,n)$. In this case put $m=m_{\xi(j,n)}$ and observe that $\lambda_j(y)=(1-\lambda_{j+1}(y))>1-\delta^2/10\ge 9/10$ and thus
$$
\begin{aligned}
|z^*_m(f_k(x)-f_n(y))|&\ge\lambda_j(y)|z^*_m(z_{j,n})|-\lambda_{j+1}(y)
|z^*_m(z_{j+1,n})|-\lambda_{i+1}(x)|z^*_m(z_{i+1,k})|-|z^*_m(f(x)-f(y)-
\lambda_i(x)z_{i,k})|\\
&\ge \frac9{10}\delta -\frac{\delta^2}{10}-\frac{\delta}4-3\frac{\delta^3}{100}>\frac{\delta^3}{20}.
\end{aligned}
$$

5)  $\lambda_{j+1}(y)\le \delta^2/10$, $\xi(j,n)<\xi(i+1,k)$, $\lambda_{i+1}(x)\le \delta/4$, and $\xi(i,k)>\xi(j,n)$. In this case put $m=m_{\xi(i,k)}$ and observe that $\lambda_i(x)=1-\lambda_{i+1}(x)\ge 1-\delta/4\ge 3/4$. Then 
$$
\begin{aligned}
|z^*_m&(f_k(x)-f_n(y))|\\
&\ge\lambda_i(x)|z^*_m(z_{i,n})|-\lambda_{i+1}(x)
|z^*_m(z_{i+1,k})|-\lambda_{j+1}(y)|z^*_m(z_{j+1,n})|-|z^*_m(f(x)-f(y)-
\lambda_j(y)z_{j,k})|\\
&\ge \frac3{4}\delta -\frac{\delta}{4}-\frac{\delta^2}{10}-3\frac{\delta^3}{100}>\frac{\delta^3}{20}.
\end{aligned}
$$
\end{proof}

\section{Proof of Theorem~\ref{t2}}\label{s:pf}

Given a non-separable convex set $X$ is a Fr\'echet space $L$, consider the convex cone $$C=\{(tx,t):x\in X,\;t\in[0,+\infty)\}\subset L\times\IR$$
in $L\times\IR$ with base $X\times\{1\}$ which will be identified with $X$. 

By $\pr:C\to \IR_+$, $\pr:(x,t)\mapsto t$, we denote the projection onto the second coordinate. Observe that the map $r:C\setminus\{0\}\to X$, $r:(x,t)\mapsto x/t$, determines a retraction of $C\setminus\{0\}$ onto $X$. This retraction restricted to the set $C_{[\frac13,3]}=\pr^{-1}([\frac13,3])$ is a perfect map.

To prove that $X$ has LFAP, fix an open cover $\U$ of $X$. For each open set $U\in\U$ consider the set $\tilde U=\{(tx,t):x\in U,\; \frac13<t<3\}$. Then $\tilde \U=\{\pr^{-1}(\IR\setminus[\frac12,2]),\tilde U:U\in\U\}$ is an open cover of $C$.

By Lemma~\ref{coneSAP}, the convex cone $C$ has SAP and by Lemma~\ref{SAP-LFAP}, $C$ has LFAP. Consequently, there is a map $f:C\times \w\to C$ such that $f$ is $\tilde\U$-near to the projection $f_0:C\times\w\to C$, $f_0:(x,n)\mapsto x$, and the family  $\big(f(C\times\{n\})\big)_{n\in\w}$ is locally finite in $C$. Let $\tilde f=f|X\times\w$ and $\tilde f_0=f_0|X\times\w$. It follows from the choice of the cover $\tilde\U$ that $\tilde f(X\times\w)\subset C_{[\frac13,3]}$ and the map $g=r\circ\tilde f:X\times\w\to X$ is $\U$-near to the projection $\tilde f_0:X\times\w\to X$. 

Since the family $(\tilde f(X\times\{n\})_{n\in\w}$ is locally finite in $C_{\frac13,3]}$ and the map $r:C_{[\frac13,3]}\to X$ is perfect, the family $(r\circ \tilde f(X\times\{n\})_{n\in\w}$ is locally finite in $X$, witnessing that $X$ has LFAP.

\section{Open Problems}

The proof of Theorem~\ref{t2} heavily exploits the machinery of Banach space theory and does not work in the non-locally convex case. This leaves the following problem open:

\begin{problem} Is each non-separable completely metrizable convex AR-subset of a linear metric space homeomorphic to a Hilbert space?
\end{problem}

Even a weaker problem seems to be open:

\begin{problem} Is each complete linear metric AR-space homeomorphic to a Hilbert space?
\end{problem}  

This is true in the separable case, see \cite{DT1}, \cite{DT2}.

\end{document}